\documentclass[12pt]{article}
\usepackage{amsfonts}
\usepackage{amsmath}
\usepackage{amssymb}
\usepackage{color}
\font\smallit=cmti10  

\usepackage{amssymb,latexsym}
\usepackage{latexsym, amssymb, amsmath, amscd, amsfonts, epsfig,graphicx}
\usepackage{pifont}
\usepackage{CJK, psfrag,eepic}
\usepackage{mathrsfs}
\usepackage{txfonts}
\usepackage[perpage,symbol*]{footmisc}
\makeatletter

\renewcommand\section{\@startsection {section}{1}{\z@}
{-30pt \@plus -1ex \@minus -.2ex} {2.3ex \@plus.2ex}
{\normalfont\normalsize\bfseries}}

\renewcommand\subsection{\@startsection{subsection}{2}{\z@}
{-3.25ex\@plus -1ex \@minus -.2ex} {1.5ex \@plus .2ex}
{\normalfont\normalsize\bfseries}}

\renewcommand{\@seccntformat}[1]{\csname the#1\endcsname. }

\makeatother

\newtheorem{thm}{Theorem}[section]

\newtheorem{lemma}[thm]{Lemma}

 \DeclareMathOperator{\supp}{supp}
 
\DeclareMathOperator{\Ker}{Ker}

\def\qed{\nopagebreak \hfill $\Box$\medbreak}

\begin{document}
\title{\bf On the maximal multiplicity of long zero-sum free sequences over $C_p\oplus C_p$ }

\maketitle

\begin{center}
 \vskip 20pt
{\bf Yushuang Fan$^1$\ Linlin Wang$^1$\ Qinghai Zhong$^1$}\\
{\smallit $^1$Center for Combinatorics, LPMC-TJKLC,
Nankai University, Tianjin 300071, P.R. China}\\

\footnotetext{E-mail address:
fys850820@163.com (Y.S. Fan), \\
wanglinlin\_1986@yahoo.cn (L.L. Wang), zhongqinghai@yahoo.com.cn
(Q.H. Zhong) }
\end{center}

\begin{abstract}
In this paper, we point out that the method used in [Acta Arith.
128(2007) 245-279] can be modified slightly to obtain the following
result. Let $\varepsilon \in (0,\frac 14)$ and $c>0$, and let $p$ be
a sufficiently large prime depending  on $\varepsilon$ and $c$. Then
every zero-sumfree sequence $S$ over $C_p\oplus C_p$ of length
$|S|\geq 2p-c\sqrt{p}$ contains some element at least $\lfloor
p^{\frac14-\varepsilon}\rfloor$ times.
\end{abstract}


\section{Introduction}
The structure of long zero-sumfree sequences over a finite cyclic
group has been well studied since 1975 (See \cite{BEN},
\cite{GG2},\cite{SC}, \cite{Yuan} and \cite{GLPS}). For example, it
has been proved  by Savchev and Chen \cite{SC}, and by Yuan
\cite{Yuan} independently, that every zero-sumfree sequence over
$C_n$ of length at least $\frac n2 +1$ is a partition (up to an
integer factor co-prime to $n$) of a positive integer smaller than
$n$. But for the group $G=C_n\oplus C_n$, the structure of
zero-sumfree sequences $S$ over $G$ has been determined so far only
for the case that $S$ is of the maximal length $2n-2$. In 1969, Emde
Boas and Kruyswijk \cite{Emde} conjectured that every minimal
zero-sum sequence over $C_p\oplus C_p$ of length $2p-1$ contains
some element $p-1$ times, and in 1999, Gao and Geroldinger
\cite{GG1} conjectured that the same result holds true for any group
$C_n\oplus C_n$. It is easy to see that the above conjecture is
equivalent to that every zero-sum free sequence $S$ over $G$ of
length $2n-2$ contains some element at least $n-2$ times, which
implies a complete characterization of the structure of $S$. This
conjecture has been solved quite recently by combing two results
obtained by Reiher \cite{Rei}, and by Gao, Geroldinger and
Grynkiewicz \cite{GGG}. Reiher \cite{Rei} used about forty pages to
prove that the above conjecture is true for every prime $p$, and
Gao, Geroldinger and Grynkiewicz \cite{GGG} used more than fifty
pages to prove that the above conjecture is multiple, i.e., if it is
true for $n=k$ and $n=\ell$ then it is also true for $n=k\ell$.
Unlike the case for cyclic groups, we even can't determine the
structure of zero-sumfree sequences over $C_n\oplus C_n$ of length
$2n-3$.  In this paper we shall prove the following results by
modifying the method used in \cite{GGS}.

\begin{thm}\label{theorem1} Let $\varepsilon \in (0,\frac 14)$ and $c>0$, and let $p$ be a
sufficiently large prime depending  on $\varepsilon$ and $c$. If $S$
is a zero-sumfree sequence over $C_p\oplus C_p$ of  length $|S|\geq
2p-c\sqrt{p}$, then $S$ contains some element at least $\lfloor
p^{\frac 14-\varepsilon}\rfloor$ times.
\end{thm}

\begin{thm}\label{theorem2} Let $\varepsilon \in (0,\frac 14)$ and $c>0$, and let $p$ be a
sufficiently large prime depending  on $\varepsilon$ and $c$. Let
$S$ be a
 sequence over $C_p\oplus C_p$ of  length $|S|\geq
3p-c\sqrt{p}-1$. If $S$ contains no short zero-sum subsequence then
 $S$ contains some element at least $\lfloor p^{\frac
14-\varepsilon}\rfloor$ times.
\end{thm}

\begin{thm}\label{theorem3} Let $\varepsilon \in (0,\frac 14)$ and $c>0$, and let $p$ be a
sufficiently large prime depending  on $\varepsilon$ and $c$. Let
$S$ be a
 sequence over $C_p\oplus C_p$ of  length $|S|\geq
p^2+2p-c\sqrt{p}-1$. If $S$ contains no  zero-sum subsequence of
length $p^2$ then $S$ contains some element at least $\lfloor
p^{\frac 14-\varepsilon}\rfloor$ times.
\end{thm}

\section{Notations}

Our notation and terminology are consistent with \cite{GGS}. We
briefly gather some key notions and fix the notations concerning
sequences over finite abelian groups.  Let $\mathbb{N}$ denote the
set of positive integers, $\mathbb{P}$ the set of prime integers and
$\mathbb{N}_{0}=\mathbb{N}\cup\{0\}$.
 For any two integers $a, b\in \mathbb N$, we set $[a, b]=\{x \in \mathbb{N} : a\leq x\leq b\}$.
Throughout this paper, all abelian groups will be written
additively, and for $n, r \in \mathbb N$, we denote by $C_n$ the
cyclic group of order $n$, and denote by $C_n^r$  the direct sum of
$r$ copies of $C_n$.

Let $G$ be a finite abelian group and $\exp(G)$ its exponent. Let
$\mathcal F(G)$ be the free abelian monoid, multiplicatively
written, with basis $G$. The elements of $\mathcal F(G)$ are called
sequences over $G$. We write sequences $S\in\mathcal F(G)$ in the
form
$$S=\mathop\Pi\limits_{g\in G}g^{v_g(S)},\ \text{with}\ v_g(S)\in\mathbb N_0\
\text{for all}\ g\in G.$$
 We call $v_g(G)$ the multiplicity of $g$
in $S$, and we say that $S$ contains $g$ if $v_g(S)>0.$ Further, $S$
is called squarefree if $v_g(S)\leq 1$ for all $g\in G$. The unit
element $1\in\mathcal F(G)$ is called the empty sequence. A sequence
$S_1$ is called a subsequence of $S$ if $S_1\mid S$ in $\mathcal
F(G)$. For a subset $A$ of $G$ we denote $S_{A}=\Pi_{g\in
A}g^{v_g(S)}.$  If a sequence $S\in\mathcal F(G)$ is written in the
form $S=g_1\cdot\ldots\cdot g_l,$ we tacitly assume that
$l\in\mathbb N_0$ and $g_1,\ldots,g_l\in G$.
\par For a sequence
$$S=g_1\cdot\ldots\cdot g_l=\mathop\Pi\limits_{g\in G}g^{v_g(S)}\in\mathcal
F(G),$$ we call
\begin{itemize}
\item $|S|=l=\sum_{g\in G}v_g(G)\in\mathbb N_0$ the $length$ of $S$,
\item $h(S)=\max\{v_g(S)|g\in G\}\in[0,|S|]\}$ the $maximum$ of
the $multiplicities$ of $S$,
\item$\supp(S)=\{g\in G|v_g(S)>0\}\subset G$ the $support$ of $S$,
\item$\sigma(S)=\sum_{i=1}^lg_i=\sum_{g\in G}v_g(S)g\in G$ the $sum$ of $S$,
\item$\sum_k(S)=\{\sum_{i\in I}g_i|I\subset[1,l]\ \text{with}\
|I|=k\}$ the $set$ of $k-term\  subsums$ of $S$, for all
$k\in\mathbb N,$
\item$\sum_{\leq k}(S)=\bigcup_{j\in[1,k]}\sum_j(S),\ \sum_{\geq
k}(S)=\bigcup_{j\geq k}\sum_j(S)$,
\item$\sum(S)=\sum_{\geq 1}(S)$ the $set$ of all $subsums$ of $S$.\end{itemize}

The sequence $S$ is called
\begin{itemize}
\item a $zero-sum\  sequence$ if $\sigma(S)=0$.

\item $zero-sumfree$ if $0\not\in\sum(S).$
\end{itemize}

\section{Proof of the main results}

\begin{lemma}\label{lemma1}\cite[Lemma 2.6]{WT}  Let $G$ be prime cyclic of
order $p\in\mathbb P$ and $S$ a sequence in $\mathcal F(G).$ If
$v_0(S)=0$ and $|S|=p,$ then $\sum_{\leq h(S)}(S)=G.$
\end{lemma}

\begin{lemma}\label{lemma2} \cite{JY} Let $G$ be prime cyclic of
order $p\in\mathbb P$, $S\in\mathcal F(G)$ a squarefree sequence and
$k\in [1,|S|].$
\begin{enumerate}\renewcommand{\labelenumi}{(\theenumi)}
\item $|\sum_k(S)|\geq\min\{p,k(|S|-k)+1\}$;
\item If $k=\lfloor |S|/2\rfloor$, then
$|\sum_k(S)|\geq\min\{p,(|S|^2+3)/4\}$;
\item If $|S|=\lfloor\sqrt{4p-7}\rfloor+1$ and $k=\lfloor
|S|/2\rfloor$, then $\sum_k(S)=G.$
\end{enumerate}
\end{lemma}


\begin{lemma}\label{lemma4}\cite[Lemma 4.2]{GGS}  Let $G=C_p\oplus C_p$ with $p\in\mathbb
P,$ $(e_1,e_2)$ a basis of $G$ and
$$S=\Pi_{i=1}^l(a_ie_1+b_ie_2)\in\mathcal F(G), \mbox{ where }
a_1,b_1,\ldots,a_l,b_l\in\mathbb F_p,$$ a zero-sumfree sequence of
length $|S|=l\geq p.$ Then
$$\big |\big \{\sum_{i\in I}b_i|\emptyset\neq I\subset[1,l]\ \text{with}\ \sum_{i\in I}a_i=0\big \}\big |\geq l-p+1.$$
\end{lemma}

\begin{lemma}\label{lemma5} Let
$\varepsilon \in (0,\frac 12)$,  $c>0$ and $1<r\in\mathbb N$, and
let $p$ be a sufficiently large prime depending on $\varepsilon, c$
and $r$. Let $G=C_p^r$,  and let $S$ be a sequence over $G$ of
length $|S|\geq p$. Suppose that $|S_{g+H}|\leq\lfloor cp^{\frac
12-\varepsilon}\rfloor$ holds for all subgroups $H$ of order
$p^{r-1}$ and all $g\in G$. Then $ 0\in\sum(S)$.
\end{lemma}

{\bf Proof.} Let $p$ be a sufficiently large prime depending on
$\varepsilon, c$ and $r$. Assume to the contrary that there exists a
zero-sumfree sequence
$$S=\Pi_{i=1}^sg_i\in\mathcal F(G)\ \text{of length}\
|S|=s\geq p$$ and such that $$\big |S_{g+H}\big |\leq\lfloor
cp^{\frac 12-\varepsilon}\rfloor\ \text{for any subgroup $H$ of
order}\ p^{r-1}\ \text{and any}\ g\in G.$$

Let $\hat{G}=\text{Hom}(G,\mathbb C^\times)$ be the character group
of $G$ with complex values, $\chi_0\in\hat{G}$ the principal
character, and for any $\chi\in\hat{G}$ let
$$f(\chi)=\mathop\Pi\limits_{i=1}^s(1+\chi(g_i)).$$ Clearly, we have
$f(\chi_0)=2^s$ and $$f(\chi)=1+\sum_{g\in\sum(S)}c_g\chi(g),$$
where $c_g=|\{\emptyset\neq I\subset[1,s]|\sum_{i\in I}g_i=g\}|$.
Since $S$ is zero-sumfree, we have $0\not\in\sum(S)$ and the
Orthogonality Relations(see [\cite{GH}, Lemma 5.5.2]) imply that
$$\sum_{\chi\in\hat{G}}f(\chi)=\sum_{\chi\in\hat{G}}(1+\sum_{g\in\sum(S)}c_g\chi(g))=|\hat{G}|
+\sum_{g\in\sum(S)}c_g\sum_{\chi\in\hat{G}}\chi(g)=|G|.$$ Let
$\chi\in\hat{G}\setminus\{\chi_0\}.$ We set $M=\lfloor cp^{\frac
12-\varepsilon}\rfloor$ and
$$|S|=(2k-1)M+q\ \text{with}\ q\in[0,2M-1],$$ and continue with the
following assertion:
\par A1. $|f(\chi)|\leq 2^s\exp(-\pi^2v/(2p^2))$ with $v=2M(1^2+2^2+\cdots+(k-1)^2)+qk^2.$
\par{\bf Proof of A1.} Let $j\in[-(p-1)/2,(p-1)/2]$ and $g\in G$ with $\chi(g)=\exp(2\pi
ij/p).$ Note that for any real $x$ with $|x|<\pi/2,$ we have $\cos
x\leq\exp(-x^2/2).$ Thus
\begin{equation}\label{equation6}|1+\chi(g)|=2\cos(\frac{\pi
j}p)\leq2\exp(\frac{-\pi^2j^2}{2p^2}).\end{equation} If
$H=\Ker(\chi)$, then $|H|=p^{r-1}$ and $g+H=\chi^{-1}(\exp(2\pi
ij/p)).$ Thus applying $$\big |S_{g+H}\big |\leq M$$  there are at
most $M$ elements $h\mid S$ such that $\chi(h)=\exp(2\pi ij/p).$
Consequently, the upper bound for $f(\chi)$, obtained by repeated
application of (\ref{equation6}), is maximal if the values
$0,1,-1,\ldots,k-1,-(k-1)$ are accepted $M$ times each and the
values $k,-k$ are accepted $q$ times as images of $\chi(g)$ for
$g\in\supp(S)$. Therefore
$$|f(\chi)|\leq 2^s\exp(-\pi^2v/(2p^2)).$$ Since $|S|=s=(2k-1)M+q$,
we get $k=\frac{s-q+M}{2M}$ and hence
$$\begin{aligned}
v=&2M\sum_{j=1}^{k-1}j^2+qk^2=2M\frac{(k-1)(2k-1)k}6+qk^2
\\&=\frac{(s-q-M)(s-q+M)(s-q)+3q(s-q+M)^2}{12M^2}.\end{aligned}$$
Since $q\in[0,2M-1]$ and $q\leq s,$ it follows that
$$v=\frac{s(s^2-M^2)}{12M^2}+\frac{q(2M-q)(2M+3s-2q)}{12M^2}\geq\frac{s(s^2-M^2)}{12M^2}.$$
We deduce that (here we  need $p$ sufficiently large)
\begin{equation}\label{equation7}\exp(\frac{\pi^2v}{2p^2})\geq\exp(\frac{\pi^2s(s^2-M^2)}
{24M^2p^2})>2p^r,\end{equation} where the last inequality holds
because $s\geq p$ and $p$ is sufficiently large and then
$s^2-M^2>\frac{p^2}{2}$ and
$$\frac{\pi^2s(s^2-M^2)}{24M^2p^2}>\frac{\pi^2
p^{2\varepsilon}}{ 2\cdot 24c^2}>\ln(2p^r).$$

 Therefore it follows that
$$\begin{aligned}p^r&=|G|=\sum_{\chi\in\hat{G}}f(\chi)\geq
f(\chi_0)-\sum_{\chi\neq\chi_0}|f(\chi)|\\&\geq
2^s(1-(p^r-1)\exp(\frac{-\pi^2v}{2p^2}))>2^s(1-\frac{p^r-1}{2p^r})
>2^{s-1}>p^r,\end{aligned}$$ a contradiction.

\qed

 {\bf Proof of
Theorem \ref{theorem1}}. We may assume that $c> 8$. Let $(e_1,e_2)$
be a basis of $G$ and for $i\in[1,2]$ let $\varphi_i:
G\rightarrow\langle e_i\rangle$ denote the canonical projections.
Let $\varepsilon> 0,$ and let $p$ be sufficiently large and assume
to the contrary that there exists a zero-sumfree sequence
$$S=\Pi_{i=1}^{|S|}(a_ie_1+b_ie_2)\in\mathcal F(G),\ \text{with}\
a_1,b_1,\ldots,a_s,b_s\in[0,p-1]$$ of length $|S|=s\geq
2p-c\sqrt{p}$ and with $h(S)\leq p^{\frac 14-\varepsilon}.$ Let $T$
denote a maximal squarefree subsequence of $S$ and set
$h_0=h(\varphi_1(T)).$ After renumbering if necessary we may assume
that
$$T=\Pi_{i=1}^{|T|}(a_ie_1+b_ie_2),\ a_1=\cdots=a_{h_0}=a.$$
Now we set $$W=\Pi_{i=1}^{h_0}(ae_1+b_ie_2),\ S_1=SW^{-1}$$ and
distinguish three cases.
\par Case 1: $h_0\geq\lfloor\sqrt{4p-7}\rfloor+1$. We set $k=\lfloor\sqrt{4p-7}\rfloor+1,\ l=\lfloor
k/2\rfloor$ and $$S_2=\Pi_{i=k+1}^s(a_ie_1+b_ie_2).$$ By Lemma
\ref{lemma2}(3) we have
\begin{equation}\label{equation1}\textstyle\sum_l(\Pi_{i=1}^kb_ie_2)=\langle e_2\rangle.\end{equation}
Consider the sequence $\varphi_1(S_2)=\Pi_{i=k+1}^sa_ie_1.$ Let
$v_0(\varphi_1(S_2))=t$ and after renumbering if necessary we may
set
$$W_1=\Pi_{i=k+1}^{k+1+t}(0e_1+b_ie_2),\ W_1\mid S_2.$$

Since $W_1$ is zero-sumfree, the sequence
$\varphi_2(W_1)=\Pi_{i=k+1}^{k+1+t}b_ie_2$ is a zero-sumfree
sequence over $C_p$. It follows from  Lemma \ref{lemma2}(3)  that
$|\supp (\varphi_2(W_1))|\leq \lfloor\sqrt{4p-7}\rfloor$. By the
contrary hypothesis we have that $h(\varphi_2(W_1))=h(W_1)\leq
h(S)<p^{\frac 14}$. Therefore, $t=|\varphi_2(W_1)|\leq
h(\varphi_2(W_1))|\supp (\varphi_2(W_1))|\leq p^{\frac
14}\lfloor\sqrt{4p-7}\rfloor$. Hence,
$$|\varphi_1(S_2)|-v_0(\varphi_1(S_2))=s-k-t> 2p-c\sqrt{p}
-(\lfloor\sqrt{4p-7}\rfloor+1)-p^{\frac 14}\lfloor\sqrt{4p-7}\rfloor
\geq p.$$  Thus Lemma \ref{lemma1} implies that
$\sum(\varphi_1(S_2))=\langle e_1\rangle.$ In particular, $S_2$ has
a non-empty subsequence $S_3$ such that
$\sigma(\varphi_1(S_3))=-lae_1.$ By equation (\ref{equation1}) there
is a subset $I\subset[1,k]$ such that $\sum_{i\in
I}b_ie_2=-\sigma(\varphi_2(S_3))$ and $|I|=l.$ Therefore,
$S_3\Pi_{i\in I}(ae_1+b_ie_2)$ is a non-empty zero-sum subsequence
of $S$, a contradiction.

\par Case 2: $cp^{\frac 14}\leq h_0\leq\lfloor\sqrt{4p-7}\rfloor.$
Setting $k=\lfloor h_0/2\rfloor$ and $h_1=h(\varphi_1(S_1))$ then
Lemma \ref{lemma2}(2) implies that
\begin{equation}\label{equation2}|\textstyle\sum_k(\Pi_{i=1}^{h_0}b_ie_2)|\geq\frac{h_0^2+3}4\end{equation}
and by the assumption of Case 2 we get
$$h_1\leq h(\varphi_1(T))h(S)<h_0p^{1/4}.$$ Therefore,
$$|\varphi_1(S_1)|-v_0(\varphi_1(S_1))\geq
|S_1|-h_1>2p-cp^{1/2}-h_0-h_0p^{1/4}\geq p-1,$$ whence Lemma
\ref{lemma1} implies $\sum_{\leq h_1}(\varphi_1(S_1))=\langle
e_1\rangle.$ In particular, $S_1$ has a non-empty subsequence $S_4$
such that
\begin{equation}\label{equation3}\sigma(\varphi_1(S_4))=-kae_1,\
|S_4|\leq h_1.\end{equation} By equations (\ref{equation2}) and
(\ref{equation3}) we infer that
\begin{equation}\label{equation4}\sigma(S_4)+\textstyle\sum_k(W)\subset
\langle e_2\rangle,\ \big |\sigma(S_4)+\sum_k(W)\big
|\geq\frac{h_0^2+3}4.\end{equation} Set $S_5=S(S_4W)^{-1}.$ By Lemma
\ref{lemma4} we have
$$\big |\sum(S_5)\cap\langle e_2\rangle\big |\geq |S_5|-p+1.$$
Therefore, since $c>8$ and $p$ is sufficiently large,
$$\begin{aligned}&\big |\sigma(S_4)+\textstyle\sum_k(W)\big |+\big |\sum(S_5)\cap\langle e_2\rangle\big |\geq\frac{h_0^2+3}4+|S_5|-p+1
\\&\geq\frac{h_0^2+3}4+2p-cp^{1/2}-h_0p^{1/4}-h_0-p+1
\\&\geq \frac{h_0^2+3}4-h_0(p^{1/4}+1)-cp^{1/2}+1+p \\& = h_0(\frac 14 h_0-p^{1/4}-1)-cp^{1/2}+\frac 74+p \\& \geq
cp^{1/4}((\frac c4-1)p^{1/4}-1)-cp^{1/2}+\frac 74+p \geq
p.\end{aligned}$$ It follows from the Cauchy-Davenport theorem that
$$(\sigma(S_4)+\textstyle\sum_k(W))+(\sum(S_5)\cap\langle e_2\rangle)=\langle e_2\rangle,$$ whence $0\in
\sigma(S_4)+\sum_k(W))+(\sum(S_5)\cap\langle
e_2\rangle)\subset\sum(S)$, a contradiction.

\par Case 3: $h_0<cp^{1/4}$.  Note
that $|\supp(S)\cap(ae_1+\langle e_2\rangle)|=h_0.$ Thus we may
suppose that, for every subgroup $H\subset G$ with $|H|=p$ and every
$g\in G,$ we have
$$|S_{g+H}|\leq h_0h(S)\leq \lfloor cp^{\frac 12-\varepsilon}\rfloor,$$
 since otherwise we choose a
different basis $(e_1',e_2')$ of $G$ and are back to Case 1 or Case
2. Therefore applying Lemma \ref{lemma5} with  $r=2$ we deduce that
$S$ is not zero-sumfree, a contradiction.
 \qed

 \begin{lemma} \label{lemma6} (\cite{GG}, Theorem 6.7) Every
 sequence over $C_n\oplus C_n$ of length $3n-2$ contains a zero-sum
 subsequence of length $n$ or $2n$.
\end{lemma}

\bigskip
{\bf Proof of Theorem \ref{theorem2}}. Let $k=3p-2-|S|$. Then,
$$k\leq \lfloor c\sqrt{p} \rfloor-1<p.$$ Let $W=0^kS$. Then, $W$ is a
sequence over $C_p\oplus C_p$ of length $|W|=3p-2$. By Lemma
\ref{lemma6}, $W$ contains a zero-sum sequence $T$ of length $p$ or
$2p$. So, $T_1=T0^{-v_g(T)}$ is a nonempty zero-sum subsequence of
$S$. Since $S$ contains no short zero-sum subsequence, we infer that
$|T_1|>p$ and $|T|=2p$, and $T_1$ is minimal zero-sum. It follows
that $2p\geq |T_1|\geq 2p-k \geq 2p-\lfloor c\sqrt{p} \rfloor+1.$
Take an arbitrary element $g|T_1$. Therefore, $T_1g^{-1}$ is
zero-sum free and $h(T_1g^{-1})\geq p^{\frac 14-\varepsilon}$ by
Theorem \ref{theorem1}. \qed

\begin{lemma} \label{lemma6} (\cite{Gao1}, Theorem 2) Let $G$ be a
finite abelian group, and let $S$ be a sequence over $G$ of length
$|S|=|G|+k$ with $k\geq 1$. If $S$ contains no zero-sum subsequence
of length $|G|$, then there exist a subsequence $T|S$ of length
$|T|=k+1$ and an element $g\in G$ such that $g+T$ is zero-sum free.
\end{lemma}

\bigskip
{\bf Proof of Theorem \ref{theorem3}}. By Lemma \ref{lemma6}, there
exist a subsequence $T|S$ and an element $g\in C_p\oplus C_p$ such
that $g+T$ is zero-sum free and $|g+T|=|T|=|S|-p^2+1\geq
2p-c\sqrt{p}$. It follows from Theorem \ref{theorem1} that $h(S)\geq
h(T)=h(g+T)\geq p^{\frac 14-\varepsilon}$. \qed

\bigskip
{\bf Acknowledgments.} This work was supported by the PCSIRT Project
of the Ministry of Science and Technology, and the National Science
Foundation of China.

\end{document}